\documentclass[11pt, reqno]{amsart}
\usepackage[colorlinks=true, pdfstartview=FitV, linkcolor=blue, 
citecolor=blue]{hyperref}
\usepackage{multicol}

\usepackage[letterpaper, margin=1.3in]{geometry}

\usepackage{amssymb,amsmath,amsfonts}
\usepackage{bm}
\usepackage{enumitem}
\usepackage{array}

\allowdisplaybreaks

\newcommand{\N}{\mathbb{N}}
\newcommand{\R}{\mathbb{R}}

\newcommand{\Z}{\mathbb{Z}}
\newcommand{\Q}{\mathbb{Q}}
\newcommand{\rar}{\rightarrow}
\newcommand{\mc}{\mathcal}
\newcommand{\mr}{\mathrm}

\newtheorem*{theorem*}{Theorem}
\newtheorem{theorem}{Theorem}
\newtheorem{lemma}[theorem]{Lemma}

\newtheorem{proposition}[theorem]{Proposition}

\newcommand{\qm}{q_{\mathrm{min}}}
\newcommand{\cF}{\mc{F}}
\newcommand{\rI}{\mr{I}}
\newcommand{\sump}{\sideset{}{'}\sum}
\newcommand{\EE}{\mathrm{E}}
\newcommand{\pp}{\mathrm{p}}

\begin{document}

\title[Expected value of smallest denominator in a random interval]{Expected value of the smallest denominator\\ in a random interval of fixed radius}
\author{Huayang Chen, Alan Haynes}

\thanks{Research supported by NSF grant DMS 2001248.\\
\phantom{A..}MSC 2020: 11J83, 11B57}


\begin{abstract}
We compute the probability mass function of the random variable which returns the smallest denominator of a reduced fraction in a randomly chosen real interval of radius $\delta/2$. As an application, we prove that the expected value of the smallest denominator is asymptotic, as $\delta\rar 0$, to $(16/\pi^2)\delta^{-1/2}.$
\end{abstract}

\maketitle

\section{Introduction}\label{sec.Intro}
It has been known since work of Franel \cite{Fran1924} and Landau \cite{Land1924} in the 1920's that the Riemann hypothesis can be reformulated as a question about the distribution of reduced rationals in the unit interval. In light of this, it is perhaps not surprising that the solutions to some basic problems in numerical approximation carry with them a shadow of the Riemann zeta function. Here we address one such problem, which was communicated to us by E.~Sander and J.~D.~Meiss.

Section 3 of the paper \cite{SandMeis2020} by Meiss and Sander describes an experiment in which a small parameter $\delta>0$ is fixed, a point $x$ is chosen randomly and uniformly from $[0,1)$, and then the smallest denominator of a reduced fraction in the interval $(x-\delta/2,x+\delta/2)$  (actually they use radius $\delta$ instead of $\delta/2$, but this is not important), which we denote by
\begin{equation}\label{eqn.QMinDef}
	q_\text{min}(x)=\min\left\{q\in\N:\exists ~a/q\in\Q\cap (x-\delta/2,x+\delta/2) \right\},
\end{equation}
is computed. Various statistics of this quantity are numerically investigated, and it is mentioned that, apart from some related work of C.~Stewart \cite{Stew2013} (based on results of Kruyswijk and Meijer \cite{KruyMeij1977}), basic questions about the distribution of values taken by $q_{\text{min}}$ appear to be unanswered. In this paper we use tools from elementary analytic number theory to study this random variable. As an application, we derive the following formula for its expected value.
\begin{theorem}\label{thm.EVqmin}
	As $\delta\rar 0$, we have that
	\[\EE[q_\text{min}]=\frac{16}{\pi^2}\frac{1}{\delta^{1/2}}+O\left(\log^2\delta\right).\]
\end{theorem}
To prove this theorem, we will first work out the probability mass function (PMF) of $q_\text{min}$, which of course contains more information than just its expected value. In order to describe the PMF, let us define, for $\alpha,\beta>0$ and $t\ge 0$,
\[\Pi(\alpha,\beta;t)=\begin{cases}
	t&\text{if}\quad t\le\overline{\alpha}\\
	\overline{\alpha}&\text{if}\quad\overline{\alpha}<t\le\overline{\beta}\\
	\alpha+\beta-t&\text{if}\quad\overline{\beta}<t\le\alpha+\beta\\
	0&\text{if}\quad t>\alpha+\beta,
\end{cases}\]
where
\[\overline{\alpha}=\min\{\alpha,\beta\}\quad\text{and}\quad\overline{\beta}=\max\{\alpha,\beta\}.\]
In Section \ref{sec.PMF} we will prove the following result.
\begin{theorem}\label{thm.PMF}
	For $\delta\in (0,1)$, the probability mass function $\pp:\N\rar [0,1]$ of $\qm$ is supported on the set $\{1,2,\ldots ,\lfloor 1/\delta\rfloor+1\}$, and takes the values $\pp(1)=\delta$, and
	\[\pp(q)=\sump_{a+b=q}\Pi\left(\frac{1}{qa},\frac{1}{qb};\delta\right),\quad q\ge 2,\]
	where the primed sum denotes a summation over positive integers $a$ and $b$ which are relatively prime to $q$.
\end{theorem}
The PMF could in principle be used to answer questions about other statistics associated to $q_\text{min}$, especially if one is only interested in numerical approximations. However, one cannot hope in general to obtain nice asymptotic formulas with constants that are easily recognizable. As the reader will see in the proof of Theorem \ref{thm.EVqmin}, the simple value of $16/\pi^2$ appears at the end of the computation out of a somewhat mysterious cancellation between several `main terms' involving more complicated constants.

This paper is organized as follows: In Section \ref{sec.Intro} we explain notation and establish some basic results which will be used later. In Section \ref{sec.PMF} we give the proof of Theorem \ref{thm.PMF}, and in Section \ref{sec.ExpVal} the proof of Theorem \ref{thm.EVqmin}.

\noindent{\textbf{Acknowledgments:}} We would like to thank Evelyn Sander and James Meiss for several conversations which helped to spark our interest in this problem. We also thank the referee, who caught a mistake in an earlier version and also made a nice aesthetic suggestion which improved the presentation of the main result. The current version of the proof of Theorem \ref{thm.EVqmin}, including the formulation of Theorem \ref{thm.PMF}, was suggested to us by the referee. It is a streamlined version of our original proof (the original proof is available on arXiv).

\section{Preliminary results}\label{sec.Prelims}

\noindent{\bf Notation:} We use $O$ and $\ll$ to denote the standard big-oh and Vinogradov notations. If $a$ and $q$ are integers then we write $(a,q)$ for the $\mathrm{gcd}$ of $a$ and $q$. The same notation will be used for ordered pairs $(x,y)$ in $\R^2$, but the meaning should be clear from context. The symbols $\varphi,\mu,$ and $\zeta$ denote the Euler phi, M\"obius, and Riemann zeta functions, respectively. For $n\in\N$, $\sum_{d|n}$ denotes a sum over positive divisors $d$ of $n$. For $x\in\R,$ we write $\lfloor x\rfloor$ for the greatest integer less than or equal to $x$ and $\{x\}=x-\lfloor x\rfloor$ for the fractional part of $x$. Finally, $\log$ denotes the natural logarithm.

First we summarize basic results about Farey fractions, proofs of which can be found in \cite[Chapter III]{HardWrig2008}. For each $Q\in\N$, we define the Farey fractions of order $Q$ by
\[\mc{F}_Q=\left\{\frac{a}{q}~:~a,q\in\Z,~0\le a< q\le Q,~(a,q)=1\right\}.\]
The set $\mc{F}_Q$ is taken with its usual ordering in $[0,1)$. We also write
\[\mc{F}=\bigcup_{Q\in\N}\mc{F}_Q\]
for the set of all Farey fractions. Fractions $a/q\in\mc{F}$ are always assumed to be labeled so that $q\in\N$ and $(a,q)=1$. If $a/q\in\mc{F}$ and $q\ge 2$ then we write $a'/q'$ and $a''/q''$ for the unique elements of $\mc{F}\cup\{1/1\}$ with the property that
\begin{equation}\label{eqn.ConsecDef}
	\frac{a'}{q'}<\frac{a}{q}<\frac{a''}{q''}\quad\text{are consecutive in}~\mc{F}_q\cup\{1/1\}.
\end{equation}
It follows from basic properties of Farey fractions that $q=q'+q''$ and that
\[\frac{a}{q}-\frac{a'}{q'}=\frac{1}{qq'},\qquad \frac{a''}{q''}-\frac{a}{q}=\frac{1}{qq''},\quad\text{and}\quad \frac{a''}{q''}-\frac{a'}{q'}=\frac{1}{q'q''}.\]
It is evident from this that
\[aq'-a'q=a''q-aq''=a''q'-a'q''=1,\]
and it follows that $(q,q')=(q,q'')=(q',q'')=1.$ It also follows that, once $a/q$ is chosen, $q'$ and $q''$ are uniquely determined by the requirements that
\[q'=a^{-1}~\mathrm{mod}~q,\quad 1\le q'\le q,\quad\text{and}\quad q''=q-q'.\]
It is clear that as $a$ runs through all reduced residue classes modulo $q$, so does $q'$. This implies that, for each $Q\in\N$, the map
\[\phi_Q:\mc{F}_Q\rar\left\{(m,n)\in\N^2~:~1\le m\le n\le Q,~(m,n)=1\right\}\]
defined by
\[\phi_Q(a/q)=\begin{cases}
	(1,1)&\text{if}~q=1,\\
	(q',q)&\text{if}~q\ge 2,
\end{cases}\]
is a bijection.

We assume that the reader is familiar with basic properties of and identities relating the functions $\varphi,\mu,$ and $\zeta.$ We will also make use of the following elementary lemma.
\begin{lemma}\label{lem.Euler}
	If $x<y$ are real numbers and $f:[x,y]\rar\R$ is continuous and monotonic then
	\[\sum_{x<n\le y}f(n)=\int_x^yf(t)~\mathrm{d}t+O\left(|f(x)|+|f(y)|\right).\]
\end{lemma}
The proof of this lemma is an easy exercise, which we leave to the reader.



\section{Probability mass function of $q_{\text{min}}$}\label{sec.PMF}
Fix $\delta\in (0,1)$ and let $\qm :[0,1)\rar\N$ be defined by \eqref{eqn.QMinDef}. In this section we will compute the probability mass function $\pp:\N\rar [0,1]$ of the random variable $\qm$, assuming uniform probability on the sample space $[0,1)$. Every closed interval of diameter $\delta$ contains a reduced fraction with denominator no larger than
\[Q=Q(\delta)=\left\lfloor \frac{1}{\delta}\right\rfloor+1,\]
so by throwing out a subset of $x\in [0,1)$ consisting of finitely many points (a set of measure $0$), we can assume that the range of $\qm(x)$ is contained in $\{1,\ldots ,Q\}.$

First of all, we have that $\qm(x)=1$ if and only if $x$ lies in the set
\[\mathrm{I}_{0/1}:=[0,\delta/2)\cup (1-\delta/2,1),\]
and this gives that
\[\pp(1)=\delta.\]
If $x\in [0,1)\setminus \mathrm{I}_{0/1}$ is not one of the finitely many values of $x$ thrown out in the discussion above then the reduced fraction $a/q\in (x-\delta/2,x+\delta/2)$ with smallest positive denominator is an element of $\cF_Q$.

For each element $a/q\in\cF_Q$ with $q\ge 2$, let $\rI_{a/q}$ be the subset of $[0,1)$ defined by
\[\rI_{a/q}=\left\{x\in [0,1)~:~\qm(x)=q,~a/q\in(x-\delta/2,x+\delta/2)\right\}.\]
It is not difficult to show that if $a/q$ and $\tilde{a}/q$ are distinct reduced fractions with denominator $q$ then $\rI_{a/q}$ and $\rI_{\tilde{a}/q}$ are disjoint. This implies that each set $\rI_{a/q}$ is an interval, and that it is the set of $x\in[0,1)$ for which $a/q$ is the fraction with smallest positive denominator in the interval $(x-\delta/2,x+\delta/2)$. In particular, this also proves the following simple proposition.
\begin{proposition}
	For any $x\in [0,1)$, the open interval $(x-\delta/2,x+\delta/2)$ contains a unique rational with denominator $q_{\min{}}(x)$.
\end{proposition}
For simplicity in some of the following calculations, we write $\overline{\rI}_{a/q}$ for the closure of the interval $\rI_{a/q}$. Writing $a'/q'$ and $a''/q''$ for the fractions defined by \eqref{eqn.ConsecDef}, we have the following five possibilities:
\begin{enumerate}[itemsep=.1in,parsep=.1in, topsep=.01in]
	\item[(I)] If $(qq')^{-1}\ge \delta$ and $(qq'')^{-1}\ge \delta$ then
	\[\overline{\rI}_{a/q}=\left[\frac{a}{q}-\frac{\delta}{2},\frac{a}{q}+\frac{\delta}{2}\right]\quad\text{and}\quad |\rI_{a/q}|=\delta.\]
	\item[(IIa)] If $(qq')^{-1}\le \delta$ and $(qq'')^{-1}\ge \delta$ then
\[\overline{\rI}_{a/q}=\left[\frac{a'}{q'}+\frac{\delta}{2},\frac{a}{q}+\frac{\delta}{2}\right]\quad\text{and}\quad |\rI_{a/q}|=\frac{1}{qq'}.\]
	\item[(IIb)] If $(qq')^{-1}\ge \delta$ and $(qq'')^{-1}\le \delta$ then
\[\overline{\rI}_{a/q}=\left[\frac{a}{q}-\frac{\delta}{2},\frac{a''}{q''}-\frac{\delta}{2}\right]\quad\text{and}\quad |\rI_{a/q}|=\frac{1}{qq''}=\frac{1}{q(q-q')}.\]
	\item[(III)] If $(qq')^{-1}\le \delta,~(qq'')^{-1}\le \delta,~$ and $(q'q'')^{-1}\ge \delta$ then
\[\overline{\rI}_{a/q}=\left[\frac{a'}{q'}+\frac{\delta}{2},\frac{a''}{q''}-\frac{\delta}{2}\right]\quad\text{and}\quad |\rI_{a/q}|=\frac{1}{q'q''}-\delta=\frac{1}{q'(q-q')}-\delta.\]
	\item[(IV)] If $(qq')^{-1}\le \delta,~(qq'')^{-1}\le \delta,~$ and $(q'q'')^{-1}< \delta$ then
\[\rI_{a/q}=\emptyset\quad\text{and}\quad |\rI_{a/q}|=0.\]
\end{enumerate}
These possibilities have been arranged so that they correspond (with cases (IIa) and (IIb) taken together) precisely to the four cases in the definition of
\[\Pi\left(\frac{1}{qq'},\frac{1}{qq''};\delta\right).\]
It follows that,
\begin{align*}
\pp(q)&=\sum_{\substack{1\le q'\le q\\(q',q)=1}}\Pi\left(\frac{1}{qq'},\frac{1}{qq''};\delta\right).
\end{align*}
For $q=1$ we have that $\pp(q)=\delta$ and for $q\ge 2$ the above sum is equal to
\[\sump_{a+b=q}\Pi\left(\frac{1}{qa},\frac{1}{qb};\delta\right).\]
This therefore completes the proof of Theorem \ref{thm.PMF}.

\section{Expected value of $q_\text{min}$}\label{sec.ExpVal}

In this section we will prove Theorem \ref{thm.EVqmin}. Suppose that $\delta< 1/2$. Then $\pp(2)=\delta$, and the expected value of $\qm$ is
\begin{align*}
	\EE[\qm]&=\sum_{1\le q\le \frac{1}{\delta}+1}q\pp(q)\\
	&=3\delta+\sum_{q\ge 3}q\sump_{a+b=q}\Pi\left(\frac{1}{qa},\frac{1}{qb};\delta\right).
\end{align*}
Using the fact that $\Pi(\lambda\alpha,\lambda\beta;\lambda t)=\lambda\Pi(\alpha,\beta;t)$ for $\lambda>0$, the above sum is equal to
\begin{align}
	&3\delta+\sum_{q\ge 3}\sump_{a+b=q}\Pi\left(\frac{1}{a},\frac{1}{b};q\delta\right)\nonumber\\
	&\quad =3\delta+\sum_{q\ge 2}\sum_{\substack{a+b=q\\a\not= b}}\left(\sum_{d|a,b}\mu (d)\right)\Pi\left(\frac{1}{a},\frac{1}{b};q\delta\right)\nonumber\\
	&\quad =3\delta+\sum_{d\ge 1}\mu (d)\sum_{a\not= b}\Pi\left(\frac{1}{da},\frac{1}{db};d(a+b)\delta\right)\nonumber\\
	&\quad =3\delta+\sum_{d\ge 1}\frac{\mu (d)}{d}S\left(\delta d^2\right),\label{eqn.EVSum1}
\end{align}
with
\[S(t)=\sum_{a\not=b}\Pi\left(\frac{1}{a},\frac{1}{b};t(a+b)\right).\]
It is not difficult to check that
\begin{equation}\label{eqn.SEst1}
	S(t)=0\text{ for }t\ge 1/2.
\end{equation}
We now aim to show that
\begin{equation}\label{eqn.SEst2}
	S(t)=\frac{8}{3}\frac{1}{\sqrt{t}}+O\left(\log (1/t)\right)\quad\text{for}\quad 0<t\le 1/2.
\end{equation}
Let $u=\max\{a,b\}$ and $v=a+b$. The map $(a,b)\mapsto (u,v)$ is a two-to-one map from the collection of pairs of positive distinct integers $(a,b)$ onto the collection of pairs of integers $(u,v)$ satisfying
\[0<u<v<2u.\]
This observation allows us to write, for $t<1,$
\[S(t)=2\left(S_1(t)+S_2(t)+S_3(t)\right),\]
where, setting $x=1/t$, the quantities $S_1(t), S_2(t),$ and $S_3(t)$ are defined by
\begin{align*}
	S_1(t)&=\frac{1}{x}\sum_{\substack{0<u<v<2u\\uv\le x}}v,\\
	S_2(t)&=\sum_{\substack{0<u<v<2u\\v(v-u)\le x<uv}}\frac{1}{u},\quad\text{and}\\
	S_3(t)&=\sum_{\substack{0<u<v<2u\\u(v-u)\le x<v(v-u)}}\left(\frac{v}{u(v-u)}-\frac{v}{x}\right).
\end{align*}
To estimate $S_1(t)$, we have that
\begin{align*}
	xS_1(t)&=\sum_{v\le \sqrt{x}}v\sum_{\frac{v}{2}<u<v}1+\sum_{\sqrt{x}<v\le\sqrt{2x}}v\sum_{\frac{v}{2}<u\le\frac{x}{v}}1\\
	&=\sum_{v\le \sqrt{x}}\left(\frac{v^2}{2}+O(v)\right)+\sum_{\sqrt{x}<v\le\sqrt{2x}}v\left(\frac{x}{v}-\frac{v}{2}+O(1)\right)\\
	&=\frac{2}{3}\left(\sqrt{2}-1\right)x^{3/2}+O(x).
\end{align*}
To estimate $S_2(t)$, in the case when $t\le 1/2$ (so that $x\ge 2$), we have (using Lemma \ref{lem.Euler}) that
\begin{align*}
	S_2(t)&=\sum_{\sqrt{x}\le v< \sqrt{2x}}\sum_{\frac{x}{v}<u<v}\frac{1}{u}+\sum_{v\ge\sqrt{2x}}\sum_{v-\frac{x}{v}\le u<v}\frac{1}{u}\\
	&=\sum_{\sqrt{x}\le v< \sqrt{2x}}\left(\log(v^2/x)+O(v/x)\right)\\
	&\qquad +\sum_{\sqrt{2x}\le v\le x}\left(-\log (1-x/v^2)+O(1/v)\right)\\
	&=\int_{\sqrt{x}}^{\sqrt{2x}}\log(s^2/x)~\mathrm{d}s+\int_{\sqrt{2x}}^x\log\left(\frac{1}{1-x/s^2}\right)~\mathrm{d}s+O(1+\log x)\\
	&=\sqrt{x}\left(\int_1^{\sqrt{2}}\log(\sigma^2)~\mathrm{d}\sigma+\int_{\sqrt{2}}^\infty\log\left(\frac{\sigma^2}{\sigma^2-1}\right)~\mathrm{d}\sigma\right)+O(\log x).
\end{align*}
To estimate $S_3(t)$ we make the substitution $w=v-u$ to write it as
\begin{align*}
	S_3(t)&=\sum_{\substack{0<w<\frac{v}{2}\\w(v-w)\le x<vw}}\left(\frac{v}{w(v-w)}-\frac{v}{x}\right)\\
	&=\sum_{w\le \sqrt{x/2}}\sum_{\frac{x}{w}<v\le w+\frac{x}{w}}\left(\frac{1}{w}+\frac{1}{v-w}-\frac{v}{x}\right)\\
	&\qquad +\sum_{\sqrt{x/2}<w<\sqrt{x}}\sum_{2w<v\le w+\frac{x}{w}}\left(\frac{1}{w}+\frac{1}{v-w}-\frac{v}{x}\right).
\end{align*}
We deal with each of the double sums on the right hand side separately. For the first double sum, we have that
\begin{align*}
	&\sum_{w\le \sqrt{x/2}}\sum_{\frac{x}{w}<v\le w+\frac{x}{w}}\left(\frac{1}{w}+\frac{1}{v-w}-\frac{v}{x}\right)\\
	&\qquad =\sum_{w\le\sqrt{x/2}}\left(\log\left(\frac{1}{1-w^2/x}\right)-\frac{w^2}{2x}+O\left(\frac{1}{w}\right)\right)\\
	&\qquad =\int_0^{\sqrt{x/2}}\log\left(\frac{1}{1-s^2/x}\right)~\mathrm{d}s+O(1)\\
	&\qquad\qquad-\frac{1}{2x}\left(\frac{(x/2)^{3/2}}{3}+O(x)\right)+O(\log x)\\
	&\qquad =\sqrt{x}\left(\int_0^{1/\sqrt{2}}\log\left(\frac{1}{1-\sigma^2}\right)~\mathrm{d}\sigma-\frac{\sqrt{2}}{24}\right)+O(\log x).
\end{align*}
For the second double sum, we have that
\begin{align*}
	&\sum_{\sqrt{x/2}<w<\sqrt{x}}\sum_{2w<v\le w+\frac{x}{w}}\left(\frac{1}{w}+\frac{1}{v-w}-\frac{v}{x}\right)\\
	&\qquad =\sum_{\sqrt{x/2}<w<\sqrt{x}}\left(\frac{x}{2w^2}-2+\log\left(\frac{x}{w^2}\right)+\frac{3w^2}{2x}+O\left(\frac{1}{w}\right)\right)\\
	&\qquad =\frac{x}{2}\left(\frac{\sqrt{2}-1}{\sqrt{x}}+O\left(\frac{1}{x}\right)\right)-2\left(\sqrt{x}-\sqrt{x/2}+O(1)\right)\\
	&\qquad\qquad +\int_{\sqrt{x/2}}^{\sqrt{x}}\log\left(\frac{x}{s^2}\right)~\mathrm{d}s+O(1)\\
	&\qquad\qquad+\frac{1}{2x}\left(x^{3/2}-(x/2)^{3/2}+O(x)\right)+O(\log x)\\
	&\qquad=\sqrt{x}\left(\int_{1/\sqrt{2}}^1\log\left(\frac{1}{\sigma^2}\right)~\mathrm{d}\sigma+\frac{11\sqrt{2}}{8}-2\right)+O(\log x).
\end{align*}
Combining the estimates for the two double sums, we find that
\begin{align*}
	S_3(t)=\sqrt{x}\left(\int_0^{1/\sqrt{2}}\log\left(\frac{1}{1-\sigma^2}\right)~\mathrm{d}\sigma+\int_{1/\sqrt{2}}^1\log\left(\frac{1}{\sigma^2}\right)~\mathrm{d}\sigma+\frac{4\sqrt{2}}{3}-2\right)+O(\log x).
\end{align*}
Adding our estimates for the sums $S_i(t)$, we obtain
\[S(t)=\frac{C}{\sqrt{t}}+O\left(\log (1/t)\right),\]
with
\[C=4\sqrt{2}-\frac{16}{3}+2D\]
and
\begin{align*}D&=\int_1^{\sqrt{2}}\log(\sigma^2)~\mathrm{d}\sigma+\int_{\sqrt{2}}^\infty\log\left(\frac{\sigma^2}{\sigma^2-1}\right)~\mathrm{d}\sigma\\
	&\qquad+\int_0^{1/\sqrt{2}}\log\left(\frac{1}{1-\sigma^2}\right)~\mathrm{d}\sigma+\int_{1/\sqrt{2}}^1\log\left(\frac{1}{\sigma^2}\right)~\mathrm{d}\sigma.
\end{align*}
We leave it as an exercise to verify that $D=4-2\sqrt{2}$, which establishes \eqref{eqn.SEst2}. Using \eqref{eqn.SEst1} and \eqref{eqn.SEst2} in Equation \eqref{eqn.EVSum1}, we have that
\begin{align*}
	\EE[\qm]&=3\delta+\sum_{d\ge 1}\frac{\mu (d)}{d}S\left(\delta d^2\right)\\
	&=\sum_{d\le 1/\sqrt{2\delta}}\frac{\mu (d)}{d}S\left(\delta d^2\right)+O(1)\\
	&=\sum_{d\le 1/\sqrt{2\delta}}\frac{\mu (d)}{d}\left(\frac{8}{3}\frac{1}{d\sqrt{\delta}}+O\left(\log \left(\delta d^2\right)\right)\right)+O(1)\\
	&=\frac{16}{\pi^2}\frac{1}{\delta^{1/2}}+O\left(\log^2\delta\right).
\end{align*}
This completes the proof of our main result.

\vspace{.15in}

\begin{multicols}{2}
{\footnotesize
\noindent
AH: Department of Mathematics,\\
University of Houston,\\
Houston, TX, United States.\\
haynes@math.uh.edu

\noindent HC: 1 Baylor Plaza BCM,\\
Houston, TX 77030\\ United States.\\
u246242@bcm.edu
}
\end{multicols}

\end{document}